\documentclass[a4paper]{article}

\usepackage{a4wide}
\usepackage{amsmath}
\usepackage{amssymb} 
\usepackage{bm}
\usepackage{booktabs}
\usepackage{algorithm}
\usepackage{multicol}
\usepackage{multirow}
\usepackage[T1]{fontenc}
\usepackage{lmodern, textcomp}
\usepackage{amsfonts}
\usepackage{graphicx}
\usepackage{epstopdf}
\usepackage{algorithmic}
\usepackage{cleveref}
\usepackage{amsthm}
\usepackage{url}
\usepackage{color}
\usepackage{subcaption}

\title{Transformation from Bi-CG into Bi-CR using\\a residual smoothing-like scheme\thanks{\textbf{Funding:} This study was partly supported by the Grants-in-Aid for Scientific Research Program (KAKENHI) of the Japan Society for the Promotion of Science (JSPS), grant JP24K14985.}}

\date{}

\author{Arisa Kawase\thanks{Graduate School of Integrative Science and Engineering, Tokyo City University, 1-28-1 Tamazutsumi, Setagaya-ku, Tokyo 158-8557, Japan.}
\and Kensuke Aihara\thanks{Department of Computer Science, Tokyo City University, 1-28-1 Tamazutsumi, Setagaya-ku, Tokyo, 158-8557, Japan ({\tt aiharak@tcu.ac.jp}).}}

\begin{document}

\maketitle

\begin{abstract}
Residual smoothing techniques, which produce a smooth convergence behavior of linear iterative solvers, also form connections between different methods. 
For example, minimal residual smoothing can transform the residuals of the conjugate gradient (CG) method into those of the conjugate residual (CR) method for linear systems with symmetric matrices. 
In this study, we investigate whether a similar relationship can be constructed between the Bi-CG and Bi-CR methods for nonsymmetric linear systems. 
Numerical experiments regarding the aforementioned connection demonstrate the validity of our insights. 
\end{abstract}

\


\


\section{Introduction}
We consider the Krylov subspace methods \cite{Saad2003,vanderVorst2003} for solving large sparse linear systems
\begin{align}
A\bm{x} = \bm{b},\quad A \in \mathbb{R}^{n\times n},\quad \bm{b} \in \mathbb{R}^{n} \label{Ax=b}
\end{align}
that are frequently encountered in scientific computing, where $A$ is nonsingular. 
The conjugate gradient (CG) \cite{Hestenes1952} and conjugate residual (CR) \cite{Eisenstat1983} methods are representative iterative solvers that are a type of the Krylov subspace methods for systems with symmetric (positive definite) matrices. 
The CG and CR algorithms can be derived from the residual orthogonalization and minimization conditions, respectively.
For nonsymmetric matrices, the Bi-CG \cite{Fletcher1976} and Bi-CR \cite{Sogabe2005,Sogabe2009} methods are known as the basic underlying solvers and can be considered as extensions of CG and CR, respectively, to the nonsymmetric case. 

Although residual orthogonal-type methods, such as CG, may exhibit oscillations in the residuals (with respect to a certain norm), residual smoothing techniques \cite{Schonauer1987,Weiss1996} are useful for converting the residual sequence into one with a smooth convergence behavior. 
Furthermore, residual smoothing techniques also provide connections between iterative methods. 
For example, the smoothed residuals obtained by CG with minimal residual smoothing (MRS) coincide with the CR residuals \cite{Weiss1996,Walker1995}. 
Similarly, \cite{Zhou1994} indicated that when applying a certain residual smoothing scheme to Bi-CG, the same residuals as those in the quasi-minimal residual (QMR) method \cite{Freund1991} are generated. 
These connections are useful for providing new interpretations of iterative methods and for understanding their convergence properties. 
However, to the best of our knowledge, the connections between Bi-CG and Bi-CR have not been studied via residual smoothing techniques. 
Therefore, we present a heuristic (but logical) transformation from Bi-CG into Bi-CR using a residual smoothing-like scheme. 
Specifically, this study aims to generate scalar coefficients $\eta_k \in \mathbb{R}$ such that the relationship 
\begin{align}
\bm{r}_k^{bicr} = \bm{r}_{k-1}^{bicr} + \eta_k (\bm{r}_k^{bicg} - \bm{r}_{k-1}^{bicr}) \label{bicg_bicr}
\end{align}
holds for $k=1,2,\dots$, where $\bm{r}_k^{bicg}$ and $\bm{r}_k^{bicr}$ are the $k$th Bi-CG and Bi-CR residuals, respectively, and $\bm{r}_0^{bicg} = \bm{r}_0^{bicr} = \bm{b} - A\bm{x}_0$ for an initial guess $\bm{x}_0$. 

The remainder of this paper is organized as follows. 
First, the Bi-CR method is described and interpreted considering several perspectives. 
Subsequently, we consider an alternative derivation of the Bi-CR residuals based on a residual smoothing approach that exploits the existing connections between the methods. 
The Bi-CG and Bi-CR algorithms can be obtained by formally applying the CG and CR algorithms, respectively, to an extended linear system in which the coefficient matrix has a self-adjointness in a non-standard (quasi) inner product. 
Therefore, we can obtain \eqref{bicg_bicr} by incorporating an MRS-based scheme in the quasi-inner product when deriving Bi-CG from CG. 
We also compare the convergence behavior of the resulting algorithm and Bi-CR via numerical experiments to demonstrate their equivalence. 
Finally, concluding remarks and the scope of future studies are presented.

\subsection*{Notation}

Let $H \in \mathbb{R}^{n\times n}$ be a symmetric positive definite (SPD) matrix. 
For $\bm{x}, \bm{y} \in \mathbb{R}^{n}$, the $H$-inner product is defined as $(\bm{x}, \bm{y})_H := \bm{x}^\top H\bm{y}$, and $\|\cdot \|_H$ denotes the corresponding induced norm. 
When $H$ is an identity matrix $I$, it is simplified to the standard inner product $(\bm{x}, \bm{y}) = \bm{x}^\top \bm{y}$. 
The $H$-orthogonality $\bm{x} \perp_H \bm{y}$ represents $(\bm{x}, \bm{y})_H = 0$. 

The matrix $A$ is referred to as self-adjoint in the $H$-inner product when $(A\bm{x}, \bm{y})_H = (\bm{x}, A\bm{y})_H$ holds for all $\bm{x}, \bm{y}$. 
This indicates that $HA$ is symmetric. 
The CG and CR methods can be naturally applied to $A$, which has the self-adjoint characteristics indicated above, replacing the standard inner product with the $H$-inner product. 

Meanwhile, we use the following matrix in this study: 
\begin{align*}
\hat H := 
\begin{bmatrix}
O & I \\
I & O
\end{bmatrix} \in \mathbb{R}^{2n\times 2n}, 
\end{align*}
and define a $\hat H$-quasi-inner product as 
\begin{align*}
\langle \hat{\bm{x}}, \hat{\bm{y}}\rangle_{\hat H} := \hat{\bm{x}}^\top \hat H\hat{\bm{y}}
\end{align*}
for $\hat{\bm{x}}, \hat{\bm{y}} \in \mathbb{R}^{2n}$. 
Because $\hat H$ is symmetric but not positive definite, $\langle \cdot, \cdot \rangle_{\hat H}$ is not an exact inner product, and the associated norm cannot be defined; note that $\hat{\bm{x}} \neq \bm{0}$ exists such that $\langle \hat{\bm{x}}, \hat{\bm{x}} \rangle_{\hat H} \leq 0$ is satisfied. 
Regardless, $\langle \hat{\bm{x}}, \hat{\bm{y}} \rangle_{\hat H} = \langle \hat{\bm{y}}, \hat{\bm{x}} \rangle_{\hat H}$ and $\langle \alpha \hat{\bm{x}} + \hat{\bm{y}}, \hat{\bm{z}} \rangle_{\hat H} = \alpha \langle \hat{\bm{x}}, \hat{\bm{z}} \rangle_{\hat H} + \langle \hat{\bm{y}}, \hat{\bm{z}} \rangle_{\hat H}$ hold for $\hat{\bm{x}}, \hat{\bm{y}}, \hat{\bm{z}} \in \mathbb{R}^{2n}$ and $\alpha \in \mathbb{R}$; these properties are useful for formally constructing algorithms in $\langle \cdot, \cdot \rangle_{\hat H}$. 
For convenience, $\langle \hat{\bm{x}}, \hat{\bm{y}} \rangle_{\hat H} = 0$ is expressed as $\hat{\bm{x}} \perp_{\hat H} \hat{\bm{y}}$.

\section{Bi-CR method}

Based on \cite{Sogabe2005,Sogabe2009}, the conventional derivations of the Bi-CR algorithm for solving \eqref{Ax=b} with nonsymmetric matrices are described. 
The conversion from Bi-CG to Bi-CR based on the choice of the initial shadow residual is also discussed.

\subsection{Derivation of Bi-CR using $\hat H$-quasi-inner product}

We consider the following $2n$-dimensional linear system that unifies the $A\bm{x} = \bm{b}$ and $A^{\top}\tilde{\bm{x}} = \tilde{\bm{b}}$ systems. 
\begin{align}
\hat A \hat{\bm{x}} = \hat{\bm{b}},\quad \hat A := 
\begin{bmatrix}
A & O \\
O & A^\top
\end{bmatrix},\quad 
\hat{\bm{x}} := 
\begin{bmatrix}
\bm{x} \\
\tilde{\bm{x}}
\end{bmatrix},\quad 
\hat{\bm{b}} := 
\begin{bmatrix}
\bm{b} \\
\tilde{\bm{b}}
\end{bmatrix}. \label{z1}
\end{align}
Here, $\tilde{\bm{r}}_0 := \tilde{\bm{b}} - A^\top \tilde{\bm{x}}_0 \in \mathbb{R}^n$ is the initial shadow residual and is arbitrary if $(\tilde{\bm{r}}_0, \bm{r}_0) \neq 0$. 
Because $\hat H\hat A$ is symmetric, $\hat A$ is regarded as self-adjoint in the $\hat H$-quasi-inner product. 
Subsequently, formally applying the CR algorithm (constructed in $\langle \cdot, \cdot \rangle_{\hat H}$) to \eqref{z1}, and reshaping yields the original Bi-CR method described in Algorithm~\ref{algBiCR}. 

\begin{algorithm}[t]
\caption{Original Bi-CR method.}\label{algBiCR}
\begin{algorithmic}[1]
\STATE Select an initial guess $\bm{x}_{0}$, compute $\bm{r}_{0} = \bm{b} - A\bm{x}_{0}$, and choose $\tilde{\bm{r}}_{0}$.
\STATE Set $\bm{p}_{0} = \bm{r}_{0}$ and $\tilde{\bm{p}}_{0} = \tilde{\bm{r}}_{0}$, and compute $\bm{q}_{0} = A\bm{p}_{0}$.
\FOR {$k = 0,1,\dots$, until convergence}
	\STATE $\alpha_{k} = \frac{(\tilde{\bm{r}}_{k},\, A\bm{r}_{k})}{(A^\top \tilde{\bm{p}}_{k},\, \bm{q}_{k})}$
	\STATE $\bm{x}_{k+1} = \bm{x}_{k} + \alpha_{k}\bm{p}_{k},\quad \bm{r}_{k+1} = \bm{r}_{k} - \alpha_{k}\bm{q}_{k},\quad \tilde{\bm{r}}_{k+1} = \tilde{\bm{r}}_{k} - \alpha_{k}A^\top \tilde{\bm{p}}_k$
	\STATE $\beta_{k} = \frac{(\tilde{\bm{r}}_{k+1},\, A\bm{r}_{k+1})}{(\tilde{\bm{r}}_{k},\, A\bm{r}_{k})}$
	\STATE $\bm{p}_{k+1} = \bm{r}_{k+1} + \beta_{k}\bm{p}_{k},\quad \tilde{\bm{p}}_{k+1} = \tilde{\bm{r}}_{k+1} + \beta_{k}\tilde{\bm{p}}_{k}$
	\STATE $\bm{q}_{k+1} = A\bm{r}_{k+1} + \beta_{k}\bm{q}_{k}$
\ENDFOR
\end{algorithmic}
\end{algorithm}

The aforementioned derivation was presented in \cite{Sogabe2005}, and another derivation based on preconditioning was described in \cite{Sogabe2009}. 
In the latter approach, with $\hat H$ as a preconditioner, the preconditioned CR method is applied to the following symmetrized linear system: 
\begin{align*}
\begin{bmatrix}
O & A \\
A^\top & O
\end{bmatrix}
\begin{bmatrix}
\tilde{\bm{x}} \\
\bm{x}
\end{bmatrix}
=
\begin{bmatrix}
\bm{b} \\
\tilde{\bm{b}}
\end{bmatrix}. 
\end{align*} 
The resulting algorithm is also simplified to Algorithm~\ref{algBiCR}. 
Further details can be found in \cite{Sogabe2009}; further, \cite{vanderVorst2003} can be referred to the Bi-CG case.

\subsection{Converting Bi-CG to Bi-CR via the initial shadow residual}

The simple connection between the Bi-CG and Bi-CR algorithms for \eqref{Ax=b} is discussed herein. 
Let $\bm{r}_0$ ($=\bm{b}-A\bm{x}_0$) be an initial residual and $\tilde{\bm{r}}_0$ be an initial shadow residual. 
The residuals $\bm{r}_k$, search directions $\bm{p}_k$, and their shadow counterparts $\tilde{\bm{r}}_k$ and $\tilde{\bm{p}}_k$ of Bi-CG/Bi-CR can be expressed as follows: 
\begin{align*}
\bm{r}_k = R_k(A)\bm{r}_0,\quad \tilde{\bm{r}}_k = R_k(A^{\top})\tilde{\bm{r}}_0,\quad \bm{p}_k = P_k(A)\bm{r}_0,\quad \tilde{\bm{p}}_k = P_k(A^{\top})\tilde{\bm{r}}_0, 
\end{align*}
where $R_k(\lambda)$ ($R_k(0) = 1$) and $P_k(\lambda)$ are the residual and direction polynomials, respectively. 

Using $R_k(\lambda)$ and $P_k(\lambda)$, the Bi-CG and Bi-CR coefficients $\alpha_k^{bicg}$ and $\alpha_k^{bicr}$ can be rewritten as follows:  
\begin{align*}
&\alpha_{k}^{bicg} 
= \frac{(\tilde{\bm{r}}_{k}^{bicg},\, \bm{r}_{k}^{bicg})}{(\tilde{\bm{p}}_{k}^{bicg},\, A\bm{p}_{k}^{bicg})} 
= \frac{(R_k(A^{\top})\tilde{\bm{r}}_0,\, R_k(A)\bm{r}_0)}{(P_k(A^{\top})\tilde{\bm{r}}_0,\, AP_k(A)\bm{r}_0)}, \\
&\alpha_{k}^{bicr} 
= \frac{(\tilde{\bm{r}}_{k}^{bicr},\, A\bm{r}_{k}^{bicr})}{(A^{\top}\tilde{\bm{p}}_{k}^{bicr},\, A\bm{p}_{k}^{bicr})}
= \frac{(R_k(A^{\top})\tilde{\bm{r}}_0,\, AR_k(A)\bm{r}_0)}{(A^{\top}P_k(A^{\top})\tilde{\bm{r}}_0,\, AP_k(A)\bm{r}_0)}
= \frac{(R_k(A^{\top})A^{\top}\tilde{\bm{r}}_0,\, R_k(A)\bm{r}_0)}{(P_k(A^{\top})A^{\top}\tilde{\bm{r}}_0,\, AP_k(A)\bm{r}_0)}.
\end{align*}
Thus, replacing $\tilde{\bm{r}}_0$ with $A^\top \tilde{\bm{r}}_0$ in $\alpha_k^{bicg}$ coincides with the form $\alpha_k^{bicr}$. 
A similar relationship also holds between the coefficients $\beta_k^{bicg}$ and $\beta_k^{bicr}$. 
Consequently, Bi-CG using the initial shadow residual $A^\top \tilde{\bm{r}}_0$ described in Algorithm~\ref{algBiCG} is mathematically equivalent to Bi-CR (with an initial shadow residual $\tilde{\bm{r}}_0$).

\begin{algorithm}[H]
\caption{Bi-CG method with $A^\top \tilde{\bm{r}}_0$. ($=$ Bi-CR method with $\tilde{\bm{r}}_0$)}\label{algBiCG}
\begin{algorithmic}[1]
\STATE Select an initial guess $\bm{x}_{0}$, compute $\bm{r}_{0} = \bm{b} - A\bm{x}_{0}$, and choose $\tilde{\bm{r}}_{0}$.
\STATE Compute (and overwrite) $\tilde{\bm{r}}_{0} = A^{\top}\tilde{\bm{r}}_0$, and set $\bm{p}_{0} = \bm{r}_{0}$ and $\tilde{\bm{p}}_{0} = \tilde{\bm{r}}_{0}$.
\FOR {$k = 0,1,\dots$, until convergence}
	\STATE $\alpha_{k} = \frac{(\tilde{\bm{r}}_{k},\, \bm{r}_{k})}{(\tilde{\bm{p}}_{k},\, A\bm{p}_{k})}$
	\STATE $\bm{x}_{k+1} = \bm{x}_{k} + \alpha_{k}\bm{p}_{k},\quad \bm{r}_{k+1} = \bm{r}_{k} - \alpha_{k}A\bm{p}_{k},\quad \tilde{\bm{r}}_{k+1} = \tilde{\bm{r}}_{k} - \alpha_{k}A^\top \tilde{\bm{p}}_{k}$
	\STATE $\beta_{k} = \frac{(\tilde{\bm{r}}_{k+1},\, \bm{r}_{k+1})}{(\tilde{\bm{r}}_{k},\, \bm{r}_{k})}$
	\STATE $\bm{p}_{k+1} = \bm{r}_{k+1} + \beta_{k}\bm{p}_{k},\quad \tilde{\bm{p}}_{k+1} = \tilde{\bm{r}}_{k+1} + \beta_{k}\tilde{\bm{p}}_{k}$
\ENDFOR
\end{algorithmic}
\end{algorithm}

\section{Alternative representation of Bi-CR}

We now consider an alternative transformation of Bi-CG into Bi-CR using a residual smoothing-like scheme.

\subsection{Residual smoothing technique}

We briefly review the conventional residual smoothing techniques \cite{Schonauer1987,Weiss1996,Zhou1994} to obtain a smooth convergence behavior herein. 
Let $\{\bm{r}_k\}$ and $\{\bm{x}_k\}$ be the primary sequences of the residuals and associated approximations, respectively, obtained using an iterative method. 
Subsequently, new sequences of residuals $\bm{s}_k$ and approximations $\bm{y}_k$ are generated in the form 
\begin{align}
\bm{s}_k = \bm{s}_{k-1} + \eta_k(\bm{r}_k - \bm{s}_{k-1}),\quad \bm{y}_k = \bm{y}_{k-1} + \eta_k(\bm{x}_k - \bm{y}_{k-1}), \label{Smo}
\end{align}
where $\bm{s}_0 := \bm{r}_0$, $\bm{y}_0 := \bm{x}_0$, and $\eta_k \in \mathbb{R}$ indicates the smoothing parameter. 
Typically, $\eta_k$ is selected such that the $H$-orthogonality $\bm{s}_k \perp_H (\bm{r}_k - \bm{s}_{k-1})$ is satisfied for an SPD matrix $H$, and is specifically given by the following: 
\begin{align*}
\eta_k = -\frac{(\bm{s}_{k-1},\, \bm{r}_{k} - \bm{s}_{k-1})_H}{(\bm{r}_{k} - \bm{s}_{k-1},\, \bm{r}_{k} - \bm{s}_{k-1})_H}. 
\end{align*}
This approach is equivalent to locally minimizing $\|\bm{s}_k\|_H$ and is referred to as minimal residual smoothing (MRS). 
Clearly, it holds that 
\begin{align*}
\|\bm{s}_k\|_H \leq \|\bm{r}_k\|_H,\quad \|\bm{s}_k\|_H \leq \|\bm{s}_{k-1}\|_H, 
\end{align*} 
thus, a monotonically decreasing sequence of the residual norms is obtained. 

Furthermore, certain residuals (and approximations) generated using different iterative methods can be connected in the form of \eqref{Smo}. 
For \eqref{Ax=b} with the symmetric matrices, the CG residual $\bm{r}_k^{cg}$ is known to transform into the CR residual $\bm{r}_k^{cr}$ using the aforementioned MRS with $H = I$; that is, 
\begin{align}
\bm{r}_{k}^{cr} = \bm{r}_{k-1}^{cr} + \eta_k(\bm{r}_k^{cg} - \bm{r}_{k-1}^{cr}),\quad \eta_k = -\frac{(\bm{r}_{k-1}^{cr},\, \bm{r}_{k}^{cg} - \bm{r}_{k-1}^{cr})}{(\bm{r}_{k}^{cg} - \bm{r}_{k-1}^{cr},\, \bm{r}_{k}^{cg} - \bm{r}_{k-1}^{cr})}. \label{cg_cr}
\end{align}
Moreover, the QMR residual $\bm{r}_k^{qmr}$ is known to connect to the Bi-CG residual $\bm{r}_k^{bicg}$ as follows: 
\begin{align}
\bm{r}_{k}^{qmr} = \bm{r}_{k-1}^{qmr} + \eta_k(\bm{r}_k^{bicg} - \bm{r}_{k-1}^{qmr}),\quad \eta_k = \frac{\tau_k^2}{\rho_k^2},\quad \frac{1}{\tau_k^2} = \frac{1}{\tau_{k-1}^2} + \frac{1}{\rho_k^2}, \label{bicg_qmr}
\end{align}
where $\rho_k^2 = (\bm{r}_k^{bicg},\, \bm{r}_k^{bicg})$ and $\tau_0^2 := \rho_0^2$. 
The QMR smoothing for a general sequence of residuals was developed based on this relationship. 
Further details regarding the derivations of \eqref{cg_cr} and \eqref{bicg_qmr} can be found in \cite{Weiss1996} and \cite{Zhou1994}, respectively.

\subsection{Bi-CR induced by a residual smoothing-like scheme}

As previously indicated, the Bi-CG and Bi-CR algorithms can be obtained by applying the CG and CR algorithms in $\langle \cdot, \cdot \rangle_{\hat H}$ to \eqref{z1}, respectively. 
Simultaneously, the CG residual is transformed into the CR residual using MRS, as described in \eqref{cg_cr}. 
Therefore, in the $\hat H$-quasi-inner product, incorporating the MRS-based scheme into the CG algorithm for \eqref{z1} may provide a connection between the Bi-CG and Bi-CR residuals. 

Here, we present an MRS-based scheme in $\langle \cdot, \cdot \rangle_{\hat H}$. 
Let $\hat{\bm{r}}_k \in \mathbb{R}^{2n}$ and $\hat{\bm{x}}_k \in \mathbb{R}^{2n}$ be the residuals and associated approximations obtained by applying CG in $\langle \cdot, \cdot \rangle_{\hat H}$ to \eqref{z1}. 
Subsequently, following \eqref{Smo}, we generate new residuals $\hat{\bm{s}}_k$ and approximations $\hat{\bm{y}}_k$ as follows: 
\begin{align}
\hat{\bm{s}}_k = \hat{\bm{s}}_{k-1} + \eta_k(\hat{\bm{r}}_k - \hat{\bm{s}}_{k-1}),\quad \hat{\bm{y}}_k = \hat{\bm{y}}_{k-1} + \eta_k(\hat{\bm{x}}_k - \hat{\bm{y}}_{k-1}), \label{Smo2}
\end{align}
where $\hat{\bm{s}}_0 := \hat{\bm{r}}_0$ and $\hat{\bm{y}}_0 := \hat{\bm{x}}_0$. 
Here, we determine the parameter $\eta_k \in \mathbb{R}$ such that $\hat{\bm{s}}_k \perp_{\hat H} (\hat{\bm{r}}_k - \hat{\bm{s}}_{k-1})$ is satisfied; namely, 
\begin{align}
\eta_k = -\frac{\langle \hat{\bm{s}}_{k-1},\, \hat{\bm{r}}_{k} - \hat{\bm{s}}_{k-1}\rangle_{\hat H}}{\langle \hat{\bm{r}}_{k} - \hat{\bm{s}}_{k-1},\, \hat{\bm{r}}_{k} - \hat{\bm{s}}_{k-1}\rangle_{\hat H}}. \label{eta_MRS2}
\end{align}
The specific algorithm of CG combined with \eqref{Smo2} and \eqref{eta_MRS2} is described in Algorithm~\ref{algCG}.

We then reshape Algorithm~\ref{algCG} using the same analogy of deriving the Bi-CG algorithm. 
Substituting the vectors 
\begin{align*}
\hat{\bm{s}}_k = 
\begin{bmatrix}
\bm{s}_k \\
\tilde{\bm{s}}_k
\end{bmatrix},\quad 
\hat{\bm{y}}_k = 
\begin{bmatrix}
\bm{y}_k \\
\tilde{\bm{y}}_k
\end{bmatrix},\quad 
\hat{\bm{r}}_k = 
\begin{bmatrix}
\bm{r}_k \\
\tilde{\bm{r}}_k
\end{bmatrix},\quad 
\hat{\bm{x}}_k = 
\begin{bmatrix}
\bm{x}_k \\
\tilde{\bm{x}}_k
\end{bmatrix}
\end{align*}
into \eqref{Smo2} and \eqref{eta_MRS2}, we obtain the divided recursion formulas in the standard inner product as follows: 
\begin{align*}
\bm{s}_k = \bm{s}_{k-1} + \eta_k(\bm{r}_k - \bm{s}_{k-1}),\quad \bm{y}_k = \bm{y}_{k-1} + \eta_k(\bm{x}_k - \bm{y}_{k-1}), \\
\tilde{\bm{s}}_k = \tilde{\bm{s}}_{k-1} + \eta_k(\tilde{\bm{r}}_k - \tilde{\bm{s}}_{k-1}),\quad \tilde{\bm{y}}_k = \tilde{\bm{y}}_{k-1} + \eta_k(\tilde{\bm{x}}_k - \tilde{\bm{y}}_{k-1}),
\end{align*}
where $\bm{s}_{0} = \bm{r}_{0}$, $\bm{y}_{0} = \bm{x}_{0}$, $\tilde{\bm{s}}_{0} = \tilde{\bm{r}}_{0}$, $\tilde{\bm{y}}_{0} = \tilde{\bm{x}}_{0}$, and 
\begin{align*}
\eta_k = -\frac{(\tilde{\bm{s}}_{k-1},\, \bm{r}_{k} - \bm{s}_{k-1}) + (\bm{s}_{k-1},\, \tilde{\bm{r}}_{k} - \tilde{\bm{s}}_{k-1})}{2(\tilde{\bm{r}}_{k} - \tilde{\bm{s}}_{k-1},\, \bm{r}_{k} - \bm{s}_{k-1})}. 
\end{align*}
By performing similar transformations for other vectors and scalars, the resulting algorithm is given by Algorithm~\ref{algBiCRalt}. 
Note that the auxiliary vectors $\bm{u}_k := \bm{r}_k - \bm{s}_{k-1}$ and $\tilde{\bm{u}}_k := \tilde{\bm{r}}_k - \tilde{\bm{s}}_{k-1}$ are introduced, and the update of $\tilde{\bm{y}}_k$ is omitted.

\begin{algorithm}[b]
\caption{CG method with the MRS-based scheme for \eqref{z1} in the $\hat H$-quasi-inner product.}\label{algCG}
\begin{algorithmic}[1]
\STATE Select an initial guess $\hat{\bm{x}}_{0}$, and compute $\hat{\bm{r}}_{0} = \hat{\bm{b}} -\hat{A}\hat{\bm{x}}_{0}$.
\STATE Set $\hat{\bm{p}}_{0} = \hat{\bm{r}}_{0}$, $\hat{\bm{y}}_{0} = \hat{\bm{x}}_{0}$, and $\hat{\bm{s}}_{0} = \hat{\bm{r}}_{0}$.
\FOR {$k = 0,1,\dots$, until convergence}
	\STATE $\alpha_{k} = \frac{\langle \hat{\bm{r}}_{k},\, \hat{\bm{r}}_{k}\rangle_{\hat H}}{\langle \hat{\bm{p}}_{k},\, \hat{A}\hat{\bm{p}}_{k}\rangle_{\hat H}}$
	\STATE $\hat{\bm{x}}_{k+1} = \hat{\bm{x}}_{k} + \alpha_{k}\hat{\bm{p}}_{k},\quad \hat{\bm{r}}_{k+1} = \hat{\bm{r}}_{k} - \alpha_{k}\hat{A}\hat{\bm{p}}_{k}$
	\STATE $\eta_{k+1} = - \frac{ \langle \hat{\bm{s}}_k,\, \hat{\bm{r}}_{k+1} - \hat{\bm{s}}_k \rangle_{\hat H} }{ \langle \hat{\bm{r}}_{k+1} - \hat{\bm{s}}_k,\, \hat{\bm{r}}_{k+1} - \hat{\bm{s}}_k \rangle_{\hat H} }$
	\STATE $\hat{\bm{y}}_{k+1} = \hat{\bm{y}}_k + \eta_{k+1} ( \hat{\bm{x}}_{k+1} - \hat{\bm{y}}_k ),\quad \hat{\bm{s}}_{k+1} = \hat{\bm{s}}_k + \eta_{k+1} ( \hat{\bm{r}}_{k+1} - \hat{\bm{s}}_k )$
	\STATE $\beta_{k} = \frac{\langle \hat{\bm{r}}_{k+1},\, \hat{\bm{r}}_{k+1}\rangle_{\hat H}}{\langle \hat{\bm{r}}_{k},\, \hat{\bm{r}}_{k}\rangle_{\hat H}}$
	\STATE $\hat{\bm{p}}_{k+1} = \hat{\bm{r}}_{k+1} + \beta_{k}\hat{\bm{p}}_{k}$
\ENDFOR
\end{algorithmic}
\end{algorithm}

In Algorithm~\ref{algBiCRalt}, $\bm{r}_k$ and $\bm{x}_k$ are apparently the Bi-CG residual and the associated approximation, respectively. 
Subsequently, $\bm{s}_k$ and $\bm{y}_k$ computed by a residual smoothing-like form for $\bm{r}_k$ and $\bm{x}_k$ are considered to coincide with the Bi-CR residual and associated approximation. 
In the next section, we present the numerical experiments on the model problems to demonstrate the validity of this claim.

\begin{algorithm}[t]
\caption{Alternative form of Algorithm~\ref{algCG}. (Transformation from Bi-CG into Bi-CR)}\label{algBiCRalt}
\begin{algorithmic}[1]
\STATE Select an initial guess $\bm{x}_{0}$, compute $\bm{r}_{0} = \bm{b} -A\bm{x}_{0}$, and choose $\tilde{\bm{r}}_0$.
\STATE Set $\bm{p}_{0} = \bm{r}_{0}$, $\tilde{\bm{p}}_{0} = \tilde{\bm{r}}_{0}$, $\bm{y}_{0} = \bm{x}_{0}$, $\bm{s}_{0} = \bm{r}_{0}$, and $\tilde{\bm{s}}_{0} = \tilde{\bm{r}}_{0}$.
\FOR {$k = 0,1,\dots$, until convergence}
\STATE $\alpha_{k} = \frac{(\tilde{\bm{r}}_{k},\, \bm{r}_{k})}{(\tilde{\bm{p}}_{k},\, A\bm{p}_{k})}$
\STATE $\bm{x}_{k+1} = \bm{x}_{k} + \alpha_{k}\bm{p}_{k},\quad \bm{r}_{k+1} = \bm{r}_{k} - \alpha_{k}A\bm{p}_{k},\quad \tilde{\bm{r}}_{k+1} = \tilde{\bm{r}}_{k} - \alpha_{k}A^\top \tilde{\bm{p}}_{k}$
\STATE $\bm{u}_{k+1} = \bm{r}_{k+1} - \bm{s}_k,\quad \tilde{\bm{u}}_{k+1} = \tilde{\bm{r}}_{k+1} - \tilde{\bm{s}}_k$
\STATE $\eta_{k+1} = - \frac{ ( \tilde{\bm{s}}_k,\, \bm{u}_{k+1} ) + ( \bm{s}_k,\, \tilde{\bm{u}}_{k+1} ) }{ 2( \tilde{\bm{u}}_{k+1},\, \bm{u}_{k+1} ) }$
\STATE $\bm{y}_{k+1} = \bm{y}_k + \eta_{k+1} ( \bm{x}_{k+1} - \bm{y}_k ),\quad \bm{s}_{k+1} = \bm{s}_k + \eta_{k+1} \bm{u}_{k+1},\quad \tilde{\bm{s}}_{k+1} = \tilde{\bm{s}}_k + \eta_{k+1} \tilde{\bm{u}}_{k+1}$
\STATE $\beta_{k} = \frac{(\tilde{\bm{r}}_{k+1},\, \bm{r}_{k+1})}{(\tilde{\bm{r}}_{k},\, \bm{r}_{k})}$
\STATE $\bm{p}_{k+1} = \bm{r}_{k+1} + \beta_{k}\bm{p}_{k},\quad \tilde{\bm{p}}_{k+1} = \tilde{\bm{r}}_{k+1} + \beta_{k}\tilde{\bm{p}}_{k}$
\ENDFOR
\end{algorithmic}
\end{algorithm}

\subsection{Numerical comparison}

We now apply the three variants of Bi-CR (Algorithms~\ref{algBiCR},~\ref{algBiCG}, and \ref{algBiCRalt}) to \eqref{Ax=b} and investigate whether the convergence behaviors of the residual norms coincide. 

\begin{figure}[t]
		\centering
		\includegraphics[scale=0.4]{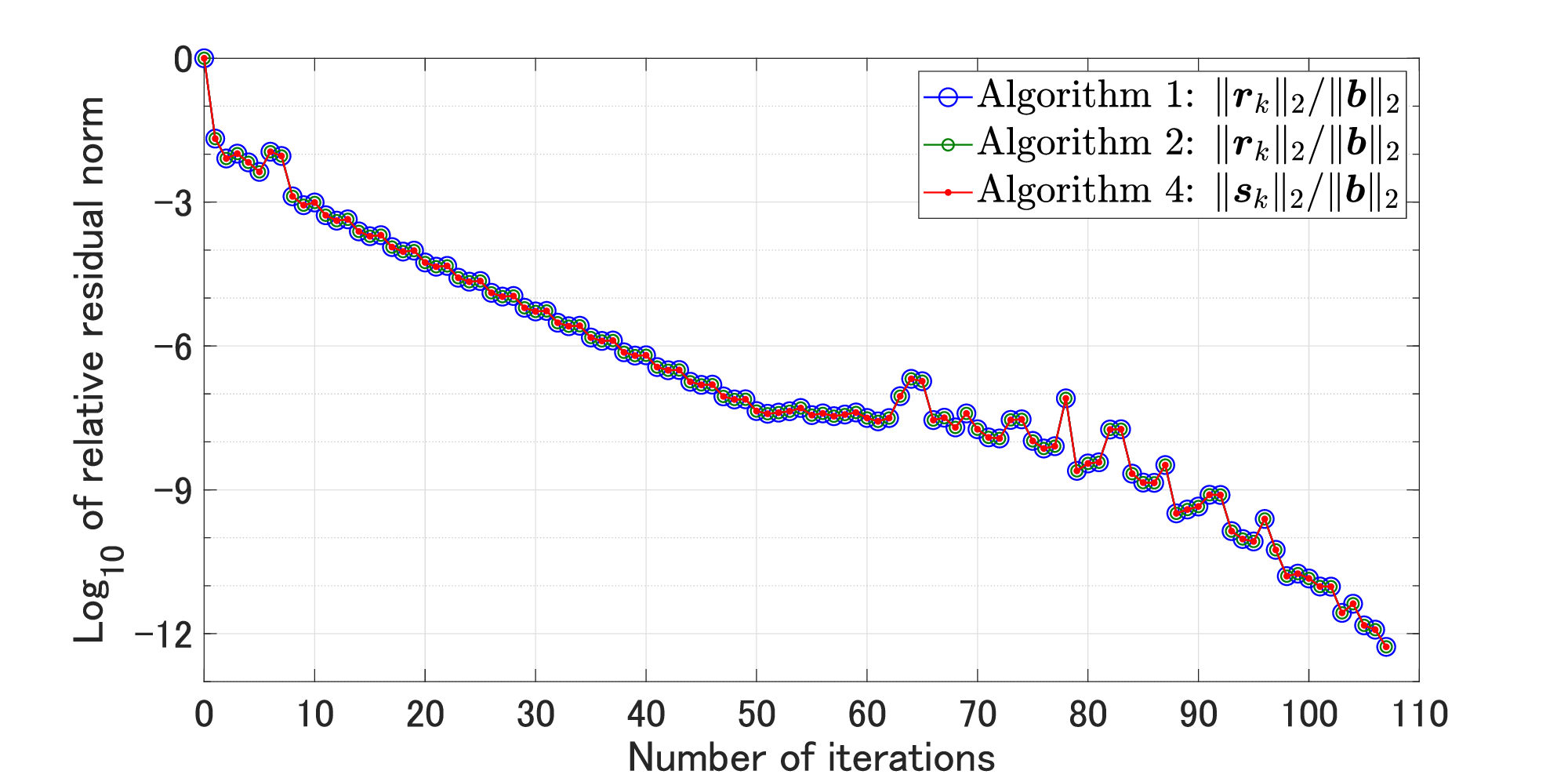}
		\caption{Convergence histories of the relative residual norms of Algorithms~\ref{algBiCR},~\ref{algBiCG}, and \ref{algBiCRalt} for the Toeplitz matrix with $\gamma = 1.2$.}\label{ex1}
		\centering
		\includegraphics[scale=0.4]{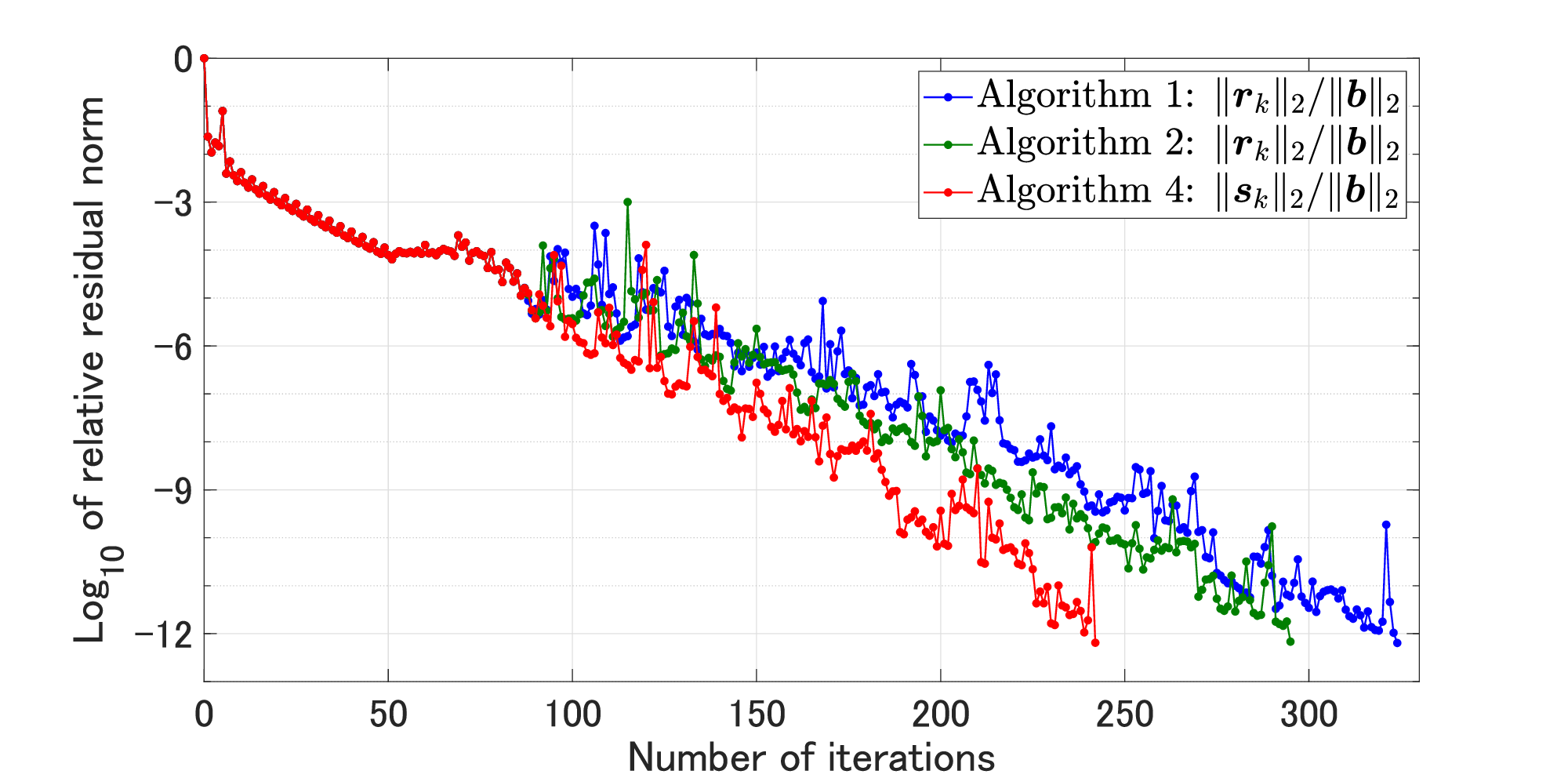}
		\caption{Convergence histories of the relative residual norms of Algorithms~\ref{algBiCR},~\ref{algBiCG}, and \ref{algBiCRalt} for the Toeplitz matrix with $\gamma = 1.5$.}\label{ex2}
\end{figure}

The numerical experiments were performed on a PC (Intel Core i7-1185G7 CPU and 32 GB of RAM) with double-precision floating-point arithmetic using MATLAB R2024a. 
Following~\cite{Sogabe2005}, as a nonsymmetric test matrix, we employed the Toeplitz matrix as follows:
\begin{align}
A := 
\begin{bmatrix}
2 & 1 & & & \\
0 & 2 & 1 &  &  \\
\gamma & 0 & 2 & 1 & \\
 & \gamma & 0 & 2 & \ddots \\
 &  & \ddots & \ddots & \ddots \\
\end{bmatrix}
\in \mathbb{R}^{200\times 200}. \label{Toep}
\end{align}
Banded Toeplitz matrices often arise from discretization of partial differential equations. 
Moreover, such matrices are also used to evaluate the convergence of the Krylov subspace methods, because their conditioning can be easily controlled. 
It has been observed that the number of Bi-CR iterations (required to converge successfully) increases for larger values of $\gamma$~\cite{Sogabe2005}. 
Thus, in the present experiments, we compare the convergence between the variants of the Bi-CR algorithms for $\gamma = 1.2$ and $\gamma = 1.5$. 
The right-hand side vector $\bm{b}$ was obtained by substituting $\bm{x}^* := [1,1,\dots,1]^\top$ into $\bm{b} = A\bm{x}^*$. 
The initial guess $\bm{x}_0$ and initial shadow residual $\tilde{\bm{r}}_0$ were set to $\bm{0}$ and $\bm{r}_0$ ($= \bm{b}$), respectively. 
The iterations were stopped when the relative residual 2-norms (i.e., $\|\bm{r}_k\|_2/\|\bm{b}\|_2$ and $\|\bm{s}_k\|_2/\|\bm{b}\|_2$ for Algorithms~\ref{algBiCR} and \ref{algBiCG} and Algorithm~\ref{algBiCRalt}, respectively) were less than $10^{-12}$. 

Figures~\ref{ex1} and \ref{ex2} demonstrate the convergence histories of Algorithms~\ref{algBiCR},~\ref{algBiCG}, and \ref{algBiCRalt} for \eqref{Toep} with $\gamma = 1.2$ and 1.5, respectively. 
The plots show the number of iterations on the horizontal axis versus $\log_{10}$ for the relative residual 2-norm (used in the stopping criterion) on the vertical axis. 

Figures~\ref{ex1} and \ref{ex2} reveal the following. 
For the Toeplitz matrix with $\gamma = 1.2$, the residual norms obtained using the three algorithms virtually coincide from the beginning to the end of the iterations. 
A similar coincidence was obtained for smaller $\gamma$ values and other test matrices that were easy to solve using Bi-CR. 
For $\gamma = 1.5$, the convergence behaviors of the three algorithms are nearly the same until approximately the 80th iteration, after which they behave differently. 
This phenomenon is likely owing to the accumulation and propagation of rounding errors in finite precision arithmetic. 
We also observed a similar coincidence at the beginning of the iterations, at least for a larger $\gamma$ and for many other test matrices. 
Because the main purpose of this study is to investigate the theoretical connection between Bi-CG and Bi-CR, we do not further discuss the difference in the effect of rounding errors on each algorithm in the present paper. 

Finally, the aforementioned results imply that the three algorithms are mathematically equivalent in exact arithmetic, and that the Bi-CG residuals can be transformed into Bi-CR residuals via \eqref{bicg_bicr} (if no breakdown occurs).

\section{Concluding remarks}

In this study, we discussed the relationships between the Bi-CG and Bi-CR methods for solving nonsymmetric linear systems and presented their novel connection via a residual smoothing-like scheme. 
In particular, we demonstrated that in the $\hat H$-quasi-inner product $\langle \cdot, \cdot \rangle_{\hat H}$, applying the MRS-based scheme to the CG algorithm for the extended system \eqref{z1} yields a residual smoothing form \eqref{bicg_bicr}, which transforms the Bi-CG residuals into Bi-CR residuals. 
This result is consistent with the existing connection between CG and CR via MRS as well as the derivation processes of Bi-CG and Bi-CR in $\langle \cdot, \cdot \rangle_{\hat H}$. 
Numerical experiments demonstrated that the Bi-CG algorithm with the residual smoothing-like scheme (i.e., Algorithm~\ref{algBiCRalt}) generates residuals that coincide with the original Bi-CR residuals (as well as Bi-CG residuals starting with an initial shadow residual acted by $A^\top$). 

However, as $\langle \cdot, \cdot \rangle_{\hat H}$ is not exactly an inner product, it is difficult to consider the orthogonality and the induced norm of the vectors. 
Therefore, unlike the relationship between CG and CR analyzed in \cite{Zhou1994}, we cannot provide direct proof of the transformation from Bi-CG into Bi-CR herein. 
A more theoretical analysis of this point can be presented in the future; for example, we consider that the orthogonality $(\tilde{\bm{s}}_i, A\bm{s}_j) = 0\ (i\neq j)$ holds and can be theoretically proven in Algorithm~\ref{algBiCRalt}, as in the original Bi-CR method. 
Moreover, in future studies, we will discuss the differences in the convergence behavior (in finite precision arithmetic) between the mathematically equivalent variants of Bi-CR including other refined implementations.


\begin{thebibliography}{12}%
\makeatletter
\providecommand \@ifxundefined [1]{%
 \@ifx{#1\undefined}
}%
\providecommand \@ifnum [1]{%
 \ifnum #1\expandafter \@firstoftwo
 \else \expandafter \@secondoftwo
 \fi
}%
\providecommand \@ifx [1]{%
 \ifx #1\expandafter \@firstoftwo
 \else \expandafter \@secondoftwo
 \fi
}%
\providecommand \natexlab [1]{#1}%
\providecommand \enquote  [1]{``#1''}%
\providecommand \bibnamefont  [1]{#1}%
\providecommand \bibfnamefont [1]{#1}%
\providecommand \citenamefont [1]{#1}%
\providecommand \href@noop [0]{\@secondoftwo}%
\providecommand \href [0]{\begingroup \@sanitize@url \@href}%
\providecommand \@href[1]{\@@startlink{#1}\@@href}%
\providecommand \@@href[1]{\endgroup#1\@@endlink}%
\providecommand \@sanitize@url [0]{\catcode `\\12\catcode `\$12\catcode
  `\&12\catcode `\#12\catcode `\^12\catcode `\_12\catcode `\%12\relax}%
\providecommand \@@startlink[1]{}%
\providecommand \@@endlink[0]{}%
\providecommand \url  [0]{\begingroup\@sanitize@url \@url }%
\providecommand \@url [1]{\endgroup\@href {#1}{\urlprefix }}%
\providecommand \urlprefix  [0]{URL }%
\providecommand \Eprint [0]{\href }%
\providecommand \doibase [0]{http://dx.doi.org/}%
\providecommand \selectlanguage [0]{\@gobble}%
\providecommand \bibinfo  [0]{\@secondoftwo}%
\providecommand \bibfield  [0]{\@secondoftwo}%
\providecommand \translation [1]{[#1]}%
\providecommand \BibitemOpen [0]{}%
\providecommand \bibitemStop [0]{}%
\providecommand \bibitemNoStop [0]{.\EOS\space}%
\providecommand \EOS [0]{\spacefactor3000\relax}%
\providecommand \BibitemShut  [1]{\csname bibitem#1\endcsname}%
\let\auto@bib@innerbib\@empty
\bibitem {Saad2003}%
  \BibitemOpen
  \bibfield  {author} {\bibinfo {author} {\bibfnamefont {Y.}~\bibnamefont
  {Saad}},\ }\href@noop {} {\emph {\bibinfo {title} {Iterative Methods for
  Sparse Linear Systems}}},\ \bibinfo {edition} {2nd}\ ed.\ (\bibinfo
  {publisher} {SIAM},\ \bibinfo {address} {Philadelphia},\ \bibinfo {year}
  {2003})\BibitemShut {NoStop}%
\bibitem {vanderVorst2003}%
  \BibitemOpen
  \bibfield  {author} {\bibinfo {author} {\bibfnamefont {H.~A.}\ \bibnamefont
  {van~der Vorst}},\ }\href@noop {} {\emph {\bibinfo {title} {Iterative Krylov
  Methods for Large Linear Systems}}}\ (\bibinfo  {publisher} {Cambridge
  University Press},\ \bibinfo {year} {2003})\BibitemShut {NoStop}%
\bibitem {Hestenes1952}%
  \BibitemOpen
  \bibfield  {author} {\bibinfo {author} {\bibfnamefont {M.~R.}\ \bibnamefont
  {Hestenes}}\ and\ \bibinfo {author} {\bibfnamefont {E.}~\bibnamefont
  {Stiefel}},\ }\bibfield  {title} {\enquote {\bibinfo {title} {Methods of
  conjugate gradients for solving linear systems},}\ }\href@noop {} {\bibfield
  {journal} {\bibinfo  {journal} {J. Research Nat. Bur. Standards}\ }\textbf
  {\bibinfo {volume} {49}},\ \bibinfo {pages} {409--436} (\bibinfo {year}
  {1952})}\BibitemShut {NoStop}%
\bibitem {Eisenstat1983}%
  \BibitemOpen
  \bibfield  {author} {\bibinfo {author} {\bibfnamefont {S.~C.}\ \bibnamefont
  {Eisenstat}}, \bibinfo {author} {\bibfnamefont {H.~C.}\ \bibnamefont
  {Elman}}, \ and\ \bibinfo {author} {\bibfnamefont {M.~H.}\ \bibnamefont
  {Schultz}},\ }\bibfield  {title} {\enquote {\bibinfo {title} {Variational
  iterative methods for nonsymmetric systems of linear equations},}\
  }\href@noop {} {\bibfield  {journal} {\bibinfo  {journal} {SIAM J. Numer.
  Anal.}\ }\textbf {\bibinfo {volume} {20}},\ \bibinfo {pages} {345--357}
  (\bibinfo {year} {1983})}\BibitemShut {NoStop}%
\bibitem {Fletcher1976}%
  \BibitemOpen
  \bibfield  {author} {\bibinfo {author} {\bibfnamefont {R.}~\bibnamefont
  {Fletcher}},\ }\bibfield  {title} {\enquote {\bibinfo {title} {Conjugate
  gradient methods for indefinite systems},}\ }\href@noop {} {\bibfield
  {journal} {\bibinfo  {journal} {Lecture Notes in Mathematics}\ }\textbf
  {\bibinfo {volume} {506}},\ \bibinfo {pages} {73--89} (\bibinfo {year}
  {1976})}\BibitemShut {NoStop}%
\bibitem {Sogabe2005}%
  \BibitemOpen
  \bibfield  {author} {\bibinfo {author} {\bibfnamefont {T.}~\bibnamefont
  {Sogabe}}, \bibinfo {author} {\bibfnamefont {M.}~\bibnamefont {Sugihara}}, \
  and\ \bibinfo {author} {\bibfnamefont {S.-L.}\ \bibnamefont {Zhang}},\
  }\bibfield  {title} {\enquote {\bibinfo {title} {An extension of the
  conjugate residual method for solving nonsymmetric linear systems},}\
  }\href@noop {} {\bibfield  {journal} {\bibinfo  {journal} {Trans. Japan Soc.
  Ind. Appl. Math.}\ }\textbf {\bibinfo {volume} {15}},\ \bibinfo {pages}
  {445--459} (\bibinfo {year} {2005})}\BibitemShut {NoStop}%
\bibitem {Sogabe2009}%
  \BibitemOpen
  \bibfield  {author} {\bibinfo {author} {\bibfnamefont {T.}~\bibnamefont
  {Sogabe}}, \bibinfo {author} {\bibfnamefont {M.}~\bibnamefont {Sugihara}}, \
  and\ \bibinfo {author} {\bibfnamefont {S.-L.}\ \bibnamefont {Zhang}},\
  }\bibfield  {title} {\enquote {\bibinfo {title} {An extension of the
  conjugate residual method to nonsymmetric linear systems},}\ }\href@noop {}
  {\bibfield  {journal} {\bibinfo  {journal} {J. Comput. Appl. Math.}\ }\textbf
  {\bibinfo {volume} {226}},\ \bibinfo {pages} {103--113} (\bibinfo {year}
  {2009})}\BibitemShut {NoStop}%
\bibitem {Schonauer1987}%
  \BibitemOpen
  \bibfield  {author} {\bibinfo {author} {\bibfnamefont {W.}~\bibnamefont
  {Sch\"{o}nauer}},\ }\href@noop {} {\emph {\bibinfo {title} {Scientific
  Computing on Vector Computers}}}\ (\bibinfo  {publisher} {Elsevier},\
  \bibinfo {address} {Amsterdam},\ \bibinfo {year} {1987})\BibitemShut
  {NoStop}%
\bibitem {Weiss1996}%
  \BibitemOpen
  \bibfield  {author} {\bibinfo {author} {\bibfnamefont {R.}~\bibnamefont
  {Weiss}},\ }\href@noop {} {\emph {\bibinfo {title} {Parameter-Free Iterative
  Linear Solvers}}}\ (\bibinfo  {publisher} {Akademie Verlag},\ \bibinfo
  {address} {Berlin},\ \bibinfo {year} {1996})\BibitemShut {NoStop}%
\bibitem {Walker1995}%
  \BibitemOpen
  \bibfield  {author} {\bibinfo {author} {\bibfnamefont {H.~F.}\ \bibnamefont
  {Walker}},\ }\bibfield  {title} {\enquote {\bibinfo {title} {Residual
  smoothing and peak/plateau behavior in {K}rylov subspace methods},}\
  }\href@noop {} {\bibfield  {journal} {\bibinfo  {journal} {Appl. Numer.
  Math.}\ }\textbf {\bibinfo {volume} {19}},\ \bibinfo {pages} {279--286}
  (\bibinfo {year} {1995})}\BibitemShut {NoStop}%
\bibitem {Zhou1994}%
  \BibitemOpen
  \bibfield  {author} {\bibinfo {author} {\bibfnamefont {L.}~\bibnamefont
  {Zhou}}\ and\ \bibinfo {author} {\bibfnamefont {H.~F.}\ \bibnamefont
  {Walker}},\ }\bibfield  {title} {\enquote {\bibinfo {title} {Residual
  smoothing techniques for iterative methods},}\ }\href@noop {} {\bibfield
  {journal} {\bibinfo  {journal} {SIAM J. Sci. Comput.}\ }\textbf {\bibinfo
  {volume} {15}},\ \bibinfo {pages} {297--312} (\bibinfo {year}
  {1994})}\BibitemShut {NoStop}%
\bibitem {Freund1991}%
  \BibitemOpen
  \bibfield  {author} {\bibinfo {author} {\bibfnamefont {R.~W.}\ \bibnamefont
  {Freund}}\ and\ \bibinfo {author} {\bibfnamefont {N.~M.}\ \bibnamefont
  {Nachtigal}},\ }\bibfield  {title} {\enquote {\bibinfo {title} {{QMR}: a
  quasi minimal residual method for non-{H}ermitian linear systems},}\
  }\href@noop {} {\bibfield  {journal} {\bibinfo  {journal} {Numer. Math.}\
  }\textbf {\bibinfo {volume} {60}},\ \bibinfo {pages} {315--339} (\bibinfo
  {year} {1991})}\BibitemShut {NoStop}%
\end{thebibliography}
\end{document}